\documentclass[a4paper,12pt]{article}
\usepackage[T1]{fontenc}
\usepackage{latexsym}
\usepackage{amssymb,amsmath,mathtools}
\usepackage{graphics}
\newtheorem{thm}{Theorem}

\usepackage[a4paper]{geometry}
\newcommand{\Z}{\mathbb{Z}}

\newcommand{\qedsymbol}{$\blacksquare$}

\title{\textbf{Interplay between the Chinese Remainder Theorem and the Lagrange Interpolation Formula}}
\author{Paul Jolissaint}
\date{}
\begin{document}

\maketitle

\section{Introduction}

Let $R$ denote a commutative, unital ring.  Then $R[t]$ denotes the ring of polynomials with coefficients in $R$, and
for $A,B\subset R$ we set 
\[
\begin{array}{rrl}
A+B & \coloneqq &\{a+b\colon a\in A, b\in B\},\\
AB &\coloneqq & \{a_1b_1+\cdots+a_nb_n\colon n\geq 1, a_j\in A, b_j\in B\}.
\end{array}
\]
 We recall that $A+B$ and $AB$ are ideals of $R$ if $A$ and $B$ are. The above sum and product generalize to finitely many ideals $A_1,\ldots,A_n$. Given an ideal $A$ and elements $x,y\in R$, we denote by $x\equiv y\pmod{A}$ the property that $x-y\in A$.
Finally, we denote by $R^*$ the (multiplicative) group of invertible elements (also called \textit{units}) of $R$, and, for $a\in R^*$, by $a^{-1}$ its inverse.

\bigskip
The purpose of the present note, which is partly inspired by \cite{Camargo}, is to discuss a somewhat surprising interplay between the following statements of the Chinese Remainder Theorem (CRT) and of the Lagrange Interpolation Formula (LIF). Note that the relationship between these results have already been observed, for instance in \cite{Brown,Sch}.
Let us state the versions that will be discussed here. The first one is \cite[Proposition 12.3.1]{IR} whose proof is reproduced in the next section.

\begin{thm}
(CRT) 
Let $R$ be a commutative, unital ring.
Suppose that $A_1,\ldots,A_n$ are pairwise coprime ideals in $R$: for all $i\not=j$, one has $A_i+A_j=R$. Set $A\coloneqq A_1\cdots A_n$. Then the natural map $\psi:R\rightarrow R/{A_1}\oplus \cdots \oplus R/{A_n}$ induces an isomorphism from $R/A$ onto $R/{A_1}\oplus \cdots \oplus R/{A_n}$. In particular, for all $y_1,\ldots,y_n\in R$, the system of congruences
\[
\left\{
 \begin{array}{ccl}
x & \equiv & y_1  \pmod{A_1}\\
 & \vdots &  \\
x & \equiv & y_n \pmod{A_n}
\end{array}
\right.
\]
admits a solution $x$, and if $x'$ is another solution, then $x\equiv x' \pmod{A}$.
\end{thm}

Let us state now the following version of LIF.

\begin{thm}
(LIF) Let $R$ be a commutative, unital ring and $(x_1,y_1),\ldots,(x_n,y_n)\in R\times R$ pairs such that, for all $i\not=j$, $x_i-x_j\in R^*$. Then there is a polynomial $p(t)\in R[t]$ such that $p(x_k)=y_k$ for every $k=1,\ldots,n$.
\end{thm}
Of course, Theorem 2 admits the following explicit, classical proof (see for instance \cite{Gau}): set 
\begin{equation}\label{LIF}
 p(t)\coloneqq \sum_{i=1}^n y_i \prod_{j\not=i}(x_i-x_j)^{-1}(t-x_j).   
\end{equation}
Then it is straightforward to check that $p(x_k)=y_k$ for every $k$. Furthermore, when $R$ is a field, the latter polynomial is the unique one such that $\deg(p(t))\leq n-1$.

\bigskip

In the rest of the note, we repeat the proof of \cite[Proposition 12.3.1]{IR} in the next section for the reader's convenience, we provide a proof of Theorem 2 based on Theorem 1 and finally we prove the classical version of CRT using a suitable version of LIF.

\section{Proof of Theorem 1}

Denote by $\psi_i:R\rightarrow R/A_i$ the natural quotient map $\psi_i(x)=x\pmod{A_i}$, and define 
\[
\psi: R\rightarrow R/{A_1}\oplus R/{A_2}\oplus \cdots \oplus R/{A_n}
\] 
by
$
\psi(x)=(\psi_1(x),\psi_2(x),\ldots,\psi_n(x))
$
for $ x\in R$.
We have to prove that $\psi$ is onto, and that $\ker(\psi)=A$.

Both statements rest on the following crucial observation: for each fixed $i$, expanding the product 
\[
\prod_{j\not=i}(A_i+A_j)=R
\]
we see that all summands are contained in $A_i$  except $\prod_{j\not=i}A_j$, which shows that 
\begin{equation}\label{eq1}
A_i+\prod_{j\not=i}A_j=R.
\end{equation}
Now, to show that $\psi$ is onto, let $y_1,\ldots,y_n\in R$. By Equation (\ref{eq1}), for every $i=1,\ldots,n$ we find $v_i\in A_i$ and $u_i\in \prod_{j\not=i}A_j$ such that $1=v_i+u_i$. This means that $v_i\equiv 0 \pmod{A_i}$ and that $u_i\equiv 1\pmod{A_i}$. Hence, setting 
$
x\coloneqq y_1u_1+\cdots +y_nu_n
$,
we get an element of $R$ such that $\psi(x)=(y_1\pmod{A_1},\ldots,y_n\pmod{A_n})$.\\
In order to prove that $\ker(\psi)=A$, we observe first that $\ker(\psi)=A_1\cap\cdots \cap A_n$, and that trivially $A\subset A_1\cap \cdots\cap A_n$. Thus we just need to prove that $A_1\cap \cdots \cap A_n \subset A$, which is done by induction on $n\geq 2$. For $n=2$, since $A_1+A_2=R$, there exist $a_i\in A_i$ such that $a_1+a_2=1$. Hence, for any $a\in A_1\cap A_2$, we have $a=aa_1+aa_2\in A_1A_2$. Assuming that the assumption is true for $n-1\geq 2$, we have 
$A_1\cap A_2\cap \cdots \cap A_n=A_1\cap (A_2A_3\cdots A_n)$, and by Equation (\ref{eq1}), $A_1+(A_2\cdots A_n)=R$ yields elements $x_1\in A_1$ and $x_2\in A_2\cdots A_n$ such that $x_1+x_2=1$, and we end the proof as in the case $n=2$ above.
\hfill
\qedsymbol

\section{Theorem 1 implies Theorem 2}

We apply Theorem 1 to the ring $R[t]$. Let $(x_1,y_1),\ldots,(x_n,y_n)$ be as in Theorem 2. For every $i$, let $A_i=(t-x_i)R[t]$ be the principal ideal in $R[t]$ generated by $t-x_i$. One has 
\[
1=(x_j-x_i)^{-1}(t-x_i)+(x_i-x_j)^{-1}(t-x_j)\in A_i+A_j
\]
for all $i\not=j$, which means that the ideals $A_i$'s satisfy the hypotheses of Theorem 1. Hence there exists an element $p(t)\in R[t]$ which satisfies the system of congruences
\[
\left\{
\begin{array}{ccl}
p(t) & \equiv & y_1 \pmod{A_1}\\
& \vdots  & \\
p(t) & \equiv & y_n \pmod{A_n}.
\end{array}
\right.
\]
This means in particular that, for every $k=1,\ldots,n$, there exists a polynomial $q_k(t)\in R[t]$ such that $p(t)-y_k=(t-x_k)q_k(t)$. Setting $t=x_k$ yields $p(x_k)-y_k=0$.
\hfill \qedsymbol

\bigskip\noindent
\textbf{Remark}
Uniqueness is hopeless in the degree of generality of Theorem 2: for instance, take  $R=\Z/91\Z$ and consider $p(t)\coloneqq t^2+t+1\in (\Z/91\Z)[t]$; it admits the four roots $9,16,74$ and $81$ which are the same as those of the polynomial
\[
q(t)\coloneqq 
(t-9)(t-16)(t-74)(t-81) \equiv t^4+2t^3+3t^2+2t+1\in (\Z/91)\Z[t].
\]

\section{A proof of the classical CRT inspired by LIF}

This section is strongly influenced by \cite{Camargo}.

\bigskip

The classical CRT states that, if $m=m_1\cdots m_n$ is a positive integer such that 
\[
(m_i,m_j)=1\quad \textrm{for\ all}\ i\not=j
\] 
and if $b_1,\ldots,b_n$ are arbitrary integers, then the system of congruences 
\[
\left\{
 \begin{array}{ccl}
x & \equiv & b_1  \pmod{m_1}\\
 & \vdots  & \\
x & \equiv & b_n \pmod{m_n}
\end{array}
\right.
\]
admits a solution $x$, and any two solutions differ by a multiple of $m$. 
We focus on the existence of $x$ here, that will follow from two observations:

(1) If $p(t)\in \Z[t]$ is such that $p(m_i)\equiv b_i \pmod{m_i}$ for every $i$, then $p(0)$ is a solution of the above system of congruences (see also \cite[Lemma 2]{Camargo}). Indeed, writing 
\[
p(t)=p(0)+a_1t+\cdots +a_kt^k=p(0)+t\underbrace{(a_1+\cdots+ a_kt^{k-1})}_{\eqqcolon q(t)}
\]
with $p(0),a_1,\ldots,a_k\in \Z$, we have, for every fixed $i$, 
\[
p(m_i)=p(0)+m_iq(m_i)\equiv p(0) \pmod{m_i};
\]
 hence $p(0)\equiv b_i\pmod{m_i}$ as claimed.

(2) Inspired by formula (\ref{LIF}), we observe that for $1\leq i\leq n$, 
\[
\prod_{j=1,j\not=i}^n(m_i-m_j)\equiv \prod_{j\not=i}^n (-m_j) \pmod{m_i}
\]
and as the class of $-m_j$ belongs to $(\Z/m_i\Z)^*$, one can find $r_i\in\Z$ such that 
\[
r_i\prod_{j=1,j\not=i}^n(m_i-m_j)\equiv 1 \pmod{m_i}.
\]
It suffices to set
\[
p(t)=\sum_{i=1}^n b_ir_i\prod_{j=1,j\not=i}^n (t-m_j)
\]
to get that $p(m_k)\equiv b_k \pmod{m_k}$ for every $k$.
\hfill \qedsymbol

\begin{flushleft}
     \begin{tabular}{l}
       Universit\'e de Neuch\^atel,\\
       Institut de Math\'emathiques,\\       
       Emile-Argand 11\\
       CH-2000 Neuch\^atel, Switzerland\\
       \small {pajolissaint@gmail.com}
     \end{tabular}
\end{flushleft}

\end{document}